\documentclass[12pt]{amsart}
\usepackage{amsmath}
\usepackage{amssymb}
\usepackage{latexsym}
\usepackage{amscd}
\usepackage{citesort}

\newdimen\AAdi%
\newbox\AAbo%
\def\AAk#1#2{\s_etbox\AAbo=\hbox{#2}\AAdi=\wd\AAbo\kern#1\AAdi{}}%
\def\AAr#1#2#3{\s_etbox\AAbo=\hbox{#2}\AAdi=\ht\AAbo\raise#1\AAdi\hbox{#3}}%
\font\tenmsb=msbm10 at 12pt
\font\sevenmsb=msbm7 at 8pt
\font\fivemsb=msbm5 at 6pt
\newfam\msbfam
\textfont\msbfam=\tenmsb
\scriptfont\msbfam=\sevenmsb
\scriptscriptfont\msbfam=\fivemsb
\def\Bbb#1{{\tenmsb\fam\msbfam#1}}
\newtheorem{thm}{Theorem}[section]
\newtheorem{lem}{Lemma}[section]
\newtheorem{cor}{Corollary}[section]
\newtheorem{rem}{Remark}[section]
\newtheorem{pro}{Proposition}[section]

\newcommand{\ba}{\begin{array}}
\newcommand{\ea}{\end{array}}

\newcommand{\Section}[2]{\setcounter{equation}{0}
\allowdisplaybreaks
\section[#1]{#2}}

\def\n{\nabla}
\def\bn{\overline\nabla}
\def\ir#1{\Bbb R^{#1}}
\def\hh#1{\Bbb H^{#1}}

\def\grs#1#2{\bold G_{#1,#2}}

\def\dt#1{\frac {d\,#1}{d\,t}}

\def\pd#1#2{\frac {\partial #1}{\partial #2}}

\def\a{\alpha}
\def\be{\beta}

\def\p#1{\partial #1}

\def\de{\delta}
\def\De{\Delta}

\def\ep{\varepsilon}
\def\G{\Gamma}
\def\g{\gamma}

\def\la{\lambda}
\def\La{\Lambda}
\def\om{\omega}

\def\th{\theta}

\begin{document}
\title
[Mean curvature Flow with convex Gauss image]
{Mean curvature Flow with convex Gauss image}

\author
[Y.L. Xin]{Y. L. Xin}
\address
{School of Mathematical Sciences, Fudan University,
Shanghai 200433, China
and
Key Laboratory of Mathematics for Nonlinear Sciences (Fudan University), Ministry of Education}
\email{ylxin@fudan.edu.cn}
\thanks{The research was partially supported by
NSFC  and SFECC}
\begin{abstract}
We study the mean curvature flow of complete space-like submanifolds in pseudo-Euclidean
space with bounded Gauss image, as well as that of complete submanifolds in Euclidean
space with convex Gauss image. By using the confinable property of the Gauss image under
the mean curvature flow we prove the long time existence results in both cases.
We also study the asymptotic behavior of these solutions when $t\to\infty$.
\end{abstract}
\renewcommand{\subjclassname}{%
  \textup{2000} Mathematics Subject Classification}
\date{}
\maketitle

\Section{Introduction}{Introduction}

There are many works on the mean curvature flow of hypersurfaces in  Riemannian manifolds
(see \cite{h1}, \cite{h2}, \cite{E-H1}, \cite{E-H2} for example) . The impressive features
of mean curvature flow for codimension one are as follows.

1. If the initial hypersurface $M_0\subset \ir{m+1}$ is uniformally convex, then the hypersurfaces
under the mean curvature contract smoothly to a single point in finite time and the shape
of the hypersurfaces becomes sphirical at the end of the contraction. If the ambient manifold is
a general Riemannian manifold, such a contraction is still working.

2. If the initial hypersurface $M_0\subset \ir{m+1}$ is an entire graph with linear growth, then
there is long time existence for the mean curvature flow and the shape of the hypersurfaces
becomes flat.

We know that J. Moser  \cite{Mo} proved that an entire minimal graph in $\ir{m+1}$ given by
$x^{m+1}=f(x^1, \cdots, x^m)$ with bounded gradient $|\n f|<c<\infty$ has to be hyperplane.
This  is closely related to the result of Ecker-Huisken  \cite{E-H1}, which reveals the second
feature of the mean curvature flow of hypersurfaces mentioned above. On the other hand,
Moser's result \cite{Mo} has been generalized to higher codimension in
\cite{H-J-W} \cite{FC}, and in author's joint work with J. Jost  \cite{J-X1}. Those are the
underline motivation of the present work.

It is natural to study the mean curvature flow of higher codimension. In recent years some works
have been done in  \cite{C-L1}, \cite{C-L2}, \cite{C-T}, \cite{smoczyk1},
\cite{smoczyk2}, \cite{smoczyk wang}, \cite{wang1} and  \cite{wang2} .
In the present paper we show the second feature in higher codimension. The terminology of
linear growth in \cite{E-H1} can be interpreted  as the image under the Gauss map of the
hypersurface lies in an open hemisphere. We investigate the mean curvature flow of submanifolds
with convex Gauss image naturally.

The target manifold of the Gauss map is the Grassmannian manifold in this situation. It is a symmetric
space of compact type. It has non-negative sectional curvature. The mean curvature flow is
closely related to its harmonic Gauss heat flow \cite{wang2}. The knowledge of the harmonic map
theory inspires us to investigate the dual situation firstly: mean curvature flow of a
space-like $m$-submanifold in pseudo-Euclidean space $\ir{m+n}_n$ with index $n$. Now,
the target manifold of the Gauss map is a pseudo-Grassmannian manifold. It has non-positive
sectional curvature. In the literature, rather few papers studied mean curvature flow
in an ambient Lorentzian manifold. Among them Ecker-Huisken \cite{E-H3} studied the mean curvature
flow of a compact space-like hypersurface in a Lorentzian manifold.

Whereas, there is a plenty of works on the Bernstein problem for complete space-like
submanifolds.  E.~Calabi raised the Bernstein problem for complete
space-like extremal hypersurfaces in Minkowski space $\ir{m+1}_1.$ He proved that such
hypersurfaces have to be hyperplanes when $m \le 4$ \cite{calabi}.
Cheng-Yau solved the problem for all $m$, in sharp contrast to the situation of Euclidean space \cite{Cheng-Yau}.

In \cite{Trei} and \cite{Cho-Tr}, H.~I.~Choi and A.~E.~Triebergs constructed many
complete space-like hypersurfaces with nonzero constant mean curvature
by prescribing boundary data at infinity for the Gauss map.

On the other hand, we proved \cite{X1}\cite{X2}  that for any complete space-like hypersurface $M$
with constant mean curvature in Minkowski space $\ir{m+1}_1,$ if the image under the
Gauss map $\gamma: M \to \hh{m}(-1)$ is bounded, then $M$ has to
be an $m-$plane.

Cheng-Yau's result was generalized to the higher codimension in 
\cite{Is}\cite{J-X2}. We proved in \cite{J-X2} a higher 
codimensional generalization of the above mentioned result in 
\cite{X1}.

In this paper we investigate the mean curvature deformation of
a complete  submanifold both in ambient pseudo-Euclidean space and Euclidean space.
The paper is divided by two part for two cases. The contents are organized as in the following table.

\tableofcontents

We will prove the following  main theorems  in this paper.
\begin{thm}
Let $F:M\to \ir{m+n}_n$ be a space-like complete $m-$submanifold which has bounded curvature
and bounded Gauss image. Then the evolution equation of mean curvature flow has long time smooth solution.
\end{thm}
\begin{rem} This theorem was announced in Geometric Analysis Meetings in Changsha (June) and San Diego
(July) in this summer.
\end{rem}

\begin{rem}
Our result in \cite{X1} has been refined by Xin-Ye \cite{X-Y}, independently by
Cao-Shen-Zhu \cite{C-S-Z}, as follows.

Let $M$ be a complete space-like hypersurface of constant mean
curvature in Minkowski space $\ir{m+1}_1.$ If the image of the
Gauss map $\gamma: M\to \hh{m}(-1)$ lies in a horoball in $\hh{m}(-1),$
then $M$ has to be a hyperplane.

This is the best possible result. It implies
that we may have  better result than Theorem 1.1 in codimension one case.
\end{rem}

\begin{thm}
Let $F:M\to \ir{m+n}$ be a  complete $m-$submanifold which has bounded curvature.
Suppose that the image under the Gauss map from $M$
into $\grs{m}{n}$ lies in a geodesic ball of radius
$R_0<\frac{\sqrt{2}}{12}\pi.$ Then the  evolution equation of mean curvature flow
has long time smooth solution.
\end{thm}

\begin{rem}
Compare  the above theorem with our Bernstein type result in \cite{J-X1}, better
results would be expected. It suffices to improve curvature estimates in \S3.2.
\end{rem}

\begin{rem}
From Theorem \ref{cp} and the proof of Theorem 1.2 we see that when the image under the Gauss map from $M$
into $\grs{m}{n}$ lies in a geodesic ball of radius $R_0<\frac{\sqrt{2}}4\pi.$
The equation of the mean curvature  is uniformly parabolic and has smooth solution
on some short time interval.
\end{rem}

We use the same idea to prove Theorem 1.1 and Theorem 1.2. Consider the image of the
Gauss map under the mean curvature flow. If an initial submanifold has convex Gauss
image, then deforming submanifolds under the mean curvature flow have the "confinable
property" (Theorem \ref{pbg} and Theorem \ref{cp}). This is an adequate higher codimensional
generalization of the "linear growth preserving property"  in \cite{E-H1}. In the case of
the ambient Euclidean space we have more technical issues, because of the nonnegative curvature
of the target manifold of the Gauss map.

We also study the asymptotic behavior of these solutions when $t\to\infty$, namely we study
the rescaled  mean curvature flow in both cases in \S 2.5 and \S 3.3,
respectively. The corresponding results as in \cite{E-H1} can be obtained similarly.

\bigskip

\Section{Space-like Submanifolds }{Space-like Submanifolds}
Let $\ir{m+n}_n$ be an $(m+n)$-dimensional pseudo-Euclidean space with the index
$n$. The indefinite metric is defined by
$$ds^2=\sum_{i=1}^m (dx^i)^2-\sum_{\a=m+1}^{m+n} (dx^{\a})^2$$

Let $F:M\to \ir{m+n}_n$ be a space-like $m$-submanifold in  $\ir{m+n}_n$ with
the second fundamental form $B$ defined by
$$B_{XY}\mathop{=}\limits^{def.}\left(\bn_X Y\right)^N$$
for $X, Y\in \G(TM)$. We denote $(\cdots)^T$ and $(\cdots)^N$ for the orthogonal projections into
the tangent bundle $TM$ and the normal bundle $NM$, respectively.
For $\nu\in\G(NM)$ we define the shape operator $A^\nu: TM\to TM$  by
$$A^\nu(V)= - (\bn_V\nu)^T.$$

Taking the trace of $B$ gives the mean curvature vector
$H$ of $M$ in $\ir{m+n}_n$ and
$$H\mathop{=}\limits^{def.} \text{trace}(B)
=B_{e_i e_i},$$
where $\{e_i\}$ is a local orthonormal frame field of $M.$ Here and in the sequel we use
the summation convention. The mean curvature vector is a cross-section of the normal bundle.

Choose a local Lorentzian frame field $\{e_i,e_\a\}$ along $M$ with dual frame field
$\{\omega_i,\omega_\a\}$, such that $e_i$ are tangent vectors to $M$.
We agree with the following range of indices
$$A,\,B,\,C,\cdots=1,\cdots,m+n;$$
$$i,\,j,\,k\cdots=1,\cdots, m;\;s, t=1,\,\cdots,n;\,\a,\,\be,\cdots=m+1,\cdots,m+n.$$

The induced Riemannian metric of $M$ is given by $ds_M^2 =\sum_i\omega_i^2$
and the induced structure equations of $M$ are
\begin{equation*}\begin{split}
   & d\omega_i = \omega_{ij}\wedge\omega_j,\qquad
                 \omega_{ij}+\omega_{ji} = 0,\cr
   &d\omega_{ij}= \omega_{ik}\wedge\omega_{kj}-\omega_{i\a}\wedge\omega_{\a j},\cr
   &\Omega_{ij} = d\omega_{ij}-\omega_{ik}\wedge\omega_{kj}
                = -\frac 12 R_{ijkl}\omega_k\wedge\omega_l.
\end{split}\end{equation*}
By Cartan's lemma we have
$$\omega_{\a i} = h_{\a ij}\omega_j.$$
\bigskip

\subsection{Bochner Type Formula}

We now derive the following Bochner type formula.
\begin{pro}
\begin{equation}\begin{split}
(\n^2 B)_{XY}
&=\n_X\n_Y H+ \left<B_{Xe_i},H\right>B_{Ye_i}-\left<B_{XY},B_{e_ie_j}\right>B_{e_ie_j}\\
&+2\left<B_{Xe_j},B_{Ye_i}\right>B_{e_ie_j}-\left<B_{Ye_i},B_{e_ie_j}\right>B_{Xe_j}
-\left<B_{Xe_i},B_{e_ie_j}\right>B_{Ye_j}.\label{LB}
\end{split}\end{equation}
\begin{equation}
\De|B|^2=2\,|\n B|^2+2\,\left<\n_i\n_jH,B_{ij}\right>
       +\,2\left<B_{ij},H\right>\left<B_{ik},B_{jk}\right>
       +\,2|R^\perp|^2-2\,\sum_{\a,\be}S_{\a\be}^2,\label{NLB}
\end{equation}
where $R^\perp$ denotes the curvature of the normal bundle and $S_{\a\be}=h_{\a ij}h_{\be ij}$.
\end{pro}
\begin{proof}
Choose a local orthonormal tangent frame field $\{e_i\}$ of $M$ near
$x\in M.$ Let $X, Y, \cdots$ be tangent vector fields and $\mu, \nu$
normal vector fields to $M$ near $x$  with
$$\n e_i|_x=\n_{e_i}X|_x=\n_{e_i}Y|_x=\cdots=\n_{e_i}\mu|_x
=\n_{e_i}\nu|_x=\cdots=0.$$
Thus,
$$\bn_XY|_x=\bn_XY|_x-\n_XY|_x=\left(\bn_XY\right)_x^N=B_{X\,Y},$$
$$\bn_X\mu|_x=\bn_X\mu|_x-\n_X\mu|_x=\left(\bn_X\mu\right)_x^T=-A^\mu(X),$$
$$\n_{XY}|_x\mathop{=}\limits^{def.}\n_X\n_Y|_x-\n_{\n_XY}|_x=\n_X\n_Y|_x.$$
By the Codazzi equations we have at $x$
\begin{equation}\begin{split}
(\n^2B)_{X\,Y} &=(\n_{e_i}\n_{e_i}B)_{X\,Y}\\
&=\n_{e_i}(\n_{e_i}B)_{X\,Y}
-(\n_{e_i}B)_{\n_{e_i}X\;Y}-(\n_{e_i}B)_{X\;\n_{e_i}Y}\\
&=\n_{e_i}(\n_XB)_{e_i\,Y}\\
&=(\n_X\n_{e_i}B)_{e_i\,Y}+(R_{X\,e_i}B)_{e_i\,Y}\\
&=\n_X(\n_{e_i}B)_{e_i\,Y}+(R_{X\,e_i}B)_{e_i\,Y}\\
&=\n_X\n_YH+(R_{X\, e_i}B)_{e_i\, Y}\\.\label{LB1}
\end{split}\end{equation}

Noting the Gauss equations and the Ricci equations
$$\left<R_{X\,Y}Z, W\right>+\left<Q_{X\,Y}^TZ, W\right>=0,$$
$$\left<R_{X\,Y}\mu, \nu\right>+\left<Q_{X\, Y}^N\mu, \nu\right>=0,$$
where
$$\left<Q_{X\,Y}^TZ, W\right>
=\left<B_{X\,W}, B_{Y\,Z}\right>-\left<B_{X\,Z}, B_{Y\, W}\right>,$$
$$\left<Q_{X\,Y}^N\mu, \nu\right>=\left<B_{Xe_i},\nu\right>\left<B_{Ye_i},\mu\right>
-\left<B_{Xe_i}\mu\right>\left<B_{Ye_i},\nu\right>.$$
Hence,
$$(R_{X\, e_i}B)_{e_i\, Y}
=-Q_{X\,e_i}^NB_{e_i\,Y}
   +B_{Q_{X\,e_i}^Te_i\,Y}
   +B_{e_i\,Q_{X\,e_i}^TY}.$$
Substituting the above equality into (\ref{LB1}) gives
\begin{equation}
(\n^2B)_{X\,Y}=\n_X\n_YH+A_0+B_0+C_0,\label{LB2}
\end{equation}
where
$$A_0=-Q_{X\,e_i}^NB_{e_i\,Y},\quad
B_0=B_{Q_{X\,e_i}^Te_i\,Y},\quad
C_0=B_{e_i\,Q_{X\,e_i}^TY}.$$
Now we calculate $A_0,\;B_0\;$ and $C_0$ in (\ref{LB2}).
\begin{equation*}\begin{split}
\left<A_0, \mu\right>
&=-\left<Q_{X\,e_i}^N B_{e_i\,Y}, \mu\right>\\
&=\left<B_{Xe_j},B_{Ye_i}\right>\left<B_{e_ie_j},\mu\right>
       -\left<B_{Ye_i},B_{e_ie_j}\right>\left<B_{Xe_j},\mu\right>,\\
\left<B_0, \mu\right>
&=\left<B_{Q_{X\,e_i}^Te_i\,Y}, \mu\right>
  =\left<A^\mu(Y), Q_{X\,e_i}^Te_i\right>\\
&=\left<B_{X\,A^\mu(Y)}, H\right>
   -\left<B_{X\, e_i}, B_{e_i\,A^\mu(Y)}\right>\\
&=\left<B_{Xe_i},H\right>\left<B_{Ye_i},\mu\right>
-\left<B_{Xe_i},B_{e_ie_j}\right>\left<B_{Ye_j},\mu\right>,\\
\left<C_0, \mu\right>
&=\left<B_{e_i\,Q_{X\,e_i}^TY}, \mu\right>
   =\left<A^\mu(e_i), Q_{X\,e_i}^TY\right>\\
&=\left<B_{e_ie_j},\mu\right>\left<Q^T_{Xe_i}Y,e_j\right>\\
&=\left<B_{e_ie_j},\mu\right>\left<B_{Xe_j},B_{Ye_i}\right>
-\left<B_{e_ie_j},\mu\right>\left<B_{XY},B_{e_ie_j}\right>.
\end{split}\end{equation*}
Hence,
\begin{multline}
\left<A_0+B_0+C_0,\mu\right>\\
= \left<B_{Xe_i},H\right>\left<B_{Ye_i},\mu\right>
-\left<B_{XY},B_{e_ie_j}\right>\left<B_{e_ie_j},\mu\right>\\
+2\left<B_{Xe_j},B_{Ye_i}\right>\left<B_{e_ie_j},\mu\right>
-\left<B_{Ye_i},B_{e_ie_j}\right>\left<B_{Xe_j},\mu\right>
-\left<B_{Xe_i},B_{e_ie_j}\right>\left<B_{Ye_j},\mu\right>.\label{LB3}
\end{multline}

Substituting (\ref{LB3}) into (\ref{LB2}) gives (\ref{LB}).

Denote

$$B_{ij}=B_{e_ie_j}=(\bar\n_{e_i}e_j)^N=h_{\a ij}e_\a,$$
where $\{e_\a\}$ is a local orthonormal frame field of the normal bundle near
$x\in M.$ It follows that
$$|B|^2=\sum_{i,j}\left<B_{ij},B_{ij}\right>=-\,\sum_{\a,i,j}h_{\a ij}^2\le 0.$$
We denote the absolute value of $|B|^2$ by $||B||^2,$ which is nonnegative. The same
notation for other time-like quantities. Then $||B||^2=\sum_{\a}S_{\a\a}$. It is easy to see
that
\begin{equation}
\sum_{\a\be}S_{\a\be}^2\ge \frac 1n(\sum_{\a}S_{\a\a})^2=\frac 1n ||B||^4\label{s}
\end{equation}
 Noting
\begin{equation*}\begin{split}
-\left<B_{kl},B_{ij}\right>\left<B_{ij},B_{kl}\right>
&=-(-h_{\a kl}h_{\a ij})(-h_{\be ij}h_{\be kl})\\
&=-h_{\a kl}h_{\a ij}h_{\be ij}h_{\be kl}=-\sum_{\a,\be} S_{\a\be}^2,
\end{split}\end{equation*}

\begin{equation*}\begin{split}
|R^{\perp}|^2&=\left<R_{e_ie_j}\nu_{\a},R_{e_ie_j}\nu_{\a}\right>\\
&=\left<Q_{e_i\,e_j}^N\nu_\a, Q_{e_i\,e_j}^N\nu_\a\right>\\
&=\left<B_{ik},Q_{ij}^N\nu_\a\right>\left<B_{jk},\nu_\a\right>
    -\left<B_{ik},\nu_\a\right>\left<B_{jk},Q_{ij}^N\nu_\a\right>\\
&=\left<B_{jk},\nu_\a\right>(\left<B_{il},B_{ik}\right>\left<B_{jl},\nu_\a\right>
         -\left<B_{il},\nu_\a\right>\left<B_{jl},B_{ik}\right>)\\
&\qquad-\left<B_{ik},\nu_\a\right>(\left<B_{il},B_{jk}\right>\left<B_{jl},\nu_\a\right>
                                    -\left<B_{il},\nu_\a\right>\left<B_{jl},B_{jk}\right>)\\
&=2\,\left<B_{il},B_{jk}\right>\left<B_{jl},B_{ik}\right>
  -2\,\left<B_{il},B_{ik}\right>\left<B_{jk},B_{jl}\right>,
\end{split}\end{equation*}
we have
$$\left<\n^2B,B\right>=\left<\n_i\n_jH,B_{ij}\right>
                       +\left<B_{ik},H\right>\left<B_{il},B_{kl}\right>
                        +|R^\perp|^2-\sum_{\a,\be}S_{\a\be}^2.$$
It gives (\ref{NLB}).
\end{proof}

\bigskip

\subsection{Evolution Equations}

We now consider the deformation of a submanifold under the mean curvature flow (abbreviated by MCF).
Namely, consider a one-parameter family $F_t=F(\cdot, t)$ of immersions $F_t:M\to \ir{m+n}_n$
with corresponding images $M_t=F_t(M)$ such that
\begin{equation}\begin{split}
\dt{}F(x, t)&=H(x, t),\qquad x\in M\\
F(x, 0)&=F(x)
\end{split}\label{mcf}
\end{equation}
is satisfied, where $H(x, t)$ is the mean curvature vector of $M_t$ at $F(x, t).$
Denote $e_i(t)=F_*e_i$ which is abbreviated to $e_i$ if there is no ambiguity.

\begin{equation}\begin{split}
\dt{g_{ij}}&=2\,\left<\dt{F_*e_i},F_*e_j\right>\\
&=2\,\left<\dt{\n_{e_i}F}, F_*e_j\right>=2\,\left<\n_{e_i}H,F_*e_j\right>\\
&=2\,(\bn_{e_i}\left<H,F_*e_j\right>-\left<H,\bn_{e_i}F_*e_j\right>)\\
&=-2\left<H,B_{ij}\right>.\label{tg}
\end{split}\end{equation}
It follows that
\begin{equation}
\dt{g^{ij}}=2\,g^{ik}g^{jl}\left<H,B_{kl}\right>\label{tg'}
\end{equation}
and
\begin{equation}
\dt{g}=-2\,|H|^2\,g.\label{tv},
\end{equation}
where $g=\det(g_{ij}).$  For a space-like submanifold the mean curvature
vector field is a normal vector field. It is time-like vector field.

(\ref{tv}) shows that if the initial submanifold is space-like, then it will remains
space-like for any $t$ under the mean curvature flow.

\begin{lem}
The second fundamental form and the mean curvature satisfy
\begin{equation}
\left(\dt{}-\De\right)||B||^2\le -\frac 2n ||B||^4\label{ENB},
\end{equation}
\begin{equation}
\left(\dt{}-\De\right)||H||^2\le -\frac 2n ||H||^4\label{ENH}.
\end{equation}
\end{lem}
\begin{proof}
For fixed $x_0\in M$, and $t_0$ choose orthonormal frame field $\{e_i\}$ of $M_{t_0}$
near $x_0$ which is normal at $x_0$, and orthonormal normal field $\{e_\a\}$. Then we evaluate at $x_0$ and $t_0$
in the following calculation.
\begin{equation}\begin{split}
\dt{h_{\a ij}}&=- \bn_{\dt{}}\left<\bn_{e_i}e_j,e_\a\right>\\
&=-\left<\bn_H \bn_{e_i}e_j,e_\a\right>-\left<\bn_{e_i}e_j,\bn_He_\a\right>\\
&=-\left<\bn_{e_i}\bn_{e_j}H,e_\a\right>-\left<B_{ij},\bn_He_\a\right>\\
&=-\left<\bn_{e_i}\left(\n_{e_i}H+\left(\bn_{e_j}H\right)^T\right),e_\a\right>
-\left<B_{ij},\bn_He_\a\right>\\
&=-\left<\n_{e_i}\n_{e_j}H,e_\a\right>
-\left<B_{ik}\left<\bn_{e_j}H,e_k\right>,e_\a\right>
-\left<B_{ij},\bn_He_\a\right>\\
&=-\left<\n_{e_i}\n_{e_j}H,e_\a\right>+\left<B_{ik},e_\a\right>\left<B_{jk},H\right>
-\left<B_{ij},\bn_He_\a\right>\\
&=-\left<\n_{e_i}\n_{e_j}H,e_\a\right>+h_{\a ik}h_{\be jk}H_\be
     -h_{\be ij}\left<\bn_He_\a,e_\be\right>.\label{TB1}
\end{split}\end{equation}

Since in a non-orthonormal frame field $g_{ij}=\left<F_*e_i. F_*e_j\right>$
is not a unit matrix  (except at $(x_0, t_0)$),
$$|B|^2= - g^{ik}g^{jl}h_{\a i j}\,h_{\a k l}.$$
We have at $(x_0, t_0)$
\begin{equation}\begin{split}
\dt{|B|^2}&=-\dt{g^{ik}}g^{jl}h_{\a i j}h_{\a k l}-\dt{g^{jl}}g^{ik}h_{\a i j}h_{\a k l}\\
        & - g^{ik}g^{jl}\dt{h_{\a i j}}h_{\a k l}-g^{ik}g^{jl}h_{\a ij}\dt{h_{\a k l}}\\
        &=-2\,\dt{g^{ik}}h_{\a i j}h_{\a k j}-2\,\dt{h_{\a i j}}h_{\a ij}.
\label{TB2}
\end{split}\end{equation}
From (\ref{TB1}) we have
\begin{equation}\begin{split}
\dt{h_{\a i j}}h_{\a ij}
&=h_{\a ij}\left(-\left<\n_{e_i}\n_{e_j}H, e_\a\right>
        +h_{\a ik}h_{\be jk}H_\be-h_{\be i j}\left<\bn_He_\a,e_\be\right>\right)\\
&=-h_{\a ij}\left<\n_{e_i}\n_{e_j}H, e_\a\right>
  +h_{\a ij}h_{\a ik}h_{\be jk}H_\be\\
  &\hskip1in-h_{\a ij}h_{\be ij}\left<\bn_He_\a,e_\be\right>\\
&=-h_{\a ij}\left<\n_{e_i}\n_{e_j}H, e_\a\right>
  +h_{\a i j}h_{\a ik}h_{\be j k}H_\be.
\label{TB3}
\end{split}\end{equation}
The third term of the second equality vanishes, because of the symmetric and anti-symmetric
properties in $\a$ and $\be$ simultaneously.
Noting (\ref{tg'}), we have
\begin{equation}
\dt{g^{ik}}h_{\a i j}h_{\a k j}
=2\,h_{\a ij}h_{\a k j}\left<H,B_{ik}\right>
=-2\,h_{\a ij}h_{\a k j}h_{\be ik}H_\be\label{TB4}
\end{equation}
Substituting (\ref{TB3}) and (\ref{TB4}) into (\ref{TB2})gives
\begin{equation}
\frac 12\dt{|B|^2}=\left<\n_i\n_jH,B_{ij}\right>
       +\left<B_{ij},H\right>\left<B_{ik},B_{jk}\right>.\label{TB}
\end{equation}
From (\ref{NLB}) and (\ref{TB}) we obtain
the evolution equation for the norm of the second fundamental form
\begin{equation}
\frac 12\left(\dt{}-\De\right)|B|^2=-|\n B|^2-|R^\perp|^2+\sum_{\a, \be}S_{\a\be}^2.
\label{EB}\end{equation}
From (\ref{s}) and (\ref{EB}) it follows that
\begin{equation}\begin{split}
\frac 12\left(\dt{}-\De\right)||B||^2&=-||\n B||^2-||R^\perp||^2-\sum_{\a, \be}S_{\a\be}^2\\
                                     &\le -\frac 1n ||B||^4.
\end{split}\end{equation}
$$|H|^2=\left<g^{ij}h_{\a ij}e_\a, g^{kl}h_{\be kl}e_\be\right>=-g^{ij}g^{kl}h_{\a ij}h_{\a kl}.$$
At $(x_0, t_0)$
\begin{equation*}\begin{split}
\dt{|H|^2}&=- 2 \dt{g^{ij}}g^{kl}h_{\a ij}h_{\a kl}-2 g^{ij}g^{kl}\dt{h_{\a ij}}h_{\a kl}\\
&= -2\dt{g^{ij}}h_{\a ij}H_\a - 2 g^{ij}\dt{h_{\a ij}}H_\a.
\end{split}\end{equation*}
Noting (\ref{tg'}) and (\ref{TB1})
\begin{equation*}\begin{split}
\dt{|H|^2}&=2 h_{\a ij}h_{\be ij}H_\a H_\be+2 \left<\n^2H, H\right>\\
&=2 S_{\a\be}H_\a H_\be+\De |H|^2- 2|\n_{e_i}H|^2.
\end{split}\end{equation*}
It follows that
\begin{equation*}
\left(\dt{}-\De\right)|H|^2=- 2 |\n H|^2+ 2 S_{\a\be}H_\a H_\be.
\end{equation*}
By using the Cauchy inequality
$$||H||^2=\sum_\a H_\a^2\le\sqrt{n}\sqrt{H_\a h_{\a ij}H_\be h_{\be ij}}=\sqrt{n}\sqrt{S_{\a\be}H_\a H_\be}$$
we obtain
\begin{equation}\begin{split}
\frac 12\left(\dt{}-\De\right)||H||^2&=-||\n H||^2-S_{\a\be}H_\a H_\be\\
                                     &\le -\frac 1n ||H||^4.
\end{split}\end{equation}
\end{proof}
\bigskip

\subsection{Maximum Principle and Curvature Estimates}

For a complete space-like submanifold $M\in\ir{m+n}_n$ with bounded curvature, the
Gauss equation implies its Ricci curvature is bounded from below. We can use the
well-known Omori-Yau maximum principle:

{\it Let $u$ be a $C^2-$function bounded from above on a complete manifold $M$
with Ricci curvature bounded from below. Then for any $\ep>0,$ there exists a sequence of points
$\{x_k\}\in M$, such that
$$\lim_{k\to \infty} u(x_k)=\sup u,$$
and when $k$ is sufficiently large
\begin{equation*}\begin{split}
|\n u|(x_k)<\ep,\\
\De u (x_k)<\ep.
\end{split}
\end{equation*}}
\begin{rem}
When the Ricci curvature of $M$ is bounded below by $-C\log(1+r^2\log^2(r+2))$,
where $r$ is the distance function from a fixed point $x_0\in M$, then the maximum
principle ia also applicable (see \cite{CX}).
\end{rem}

Now, we can use Omori-Yau maximum principle to do curvature estimate.

\begin{thm}
If $M_t$ is a smooth solution of (\ref{mcf}) in $[0, T)$. If $M_0$ has bounded curvature,
then there is  estimate
$$\sup_{M_t}||B||^2\le\sup_{M_0}||B||^2$$
\label{cur-es}
for $t\in [0, T)$.
\end{thm}
\begin{proof}
Let $t_0\in [0, T)$ be any time such that $\sup_{M_{t_0}}||B||^2$ is bounded.
Let $x_k(t_0)$ be a sequence of points on $M_{t_0}$ such that when $k\to \infty$
$$||B||^2(x_k(t_0))\to\sup_{M_0}||B||^2>0$$
and when $k$ is sufficiently large
\begin{equation*}\begin{split}
|\n ||B||^2|(x_k(t_0))<\ep,\\
\De ||B||^2 (x_k(t_0))<\ep
\end{split}
\end{equation*}
for any $\ep>0$. When $k$ is large enough from (\ref{ENB})
$$\left.\dt{||B||^2}\right|_{x_k(t_0)}\le \left.(\De ||B||^2-\frac 1n||B||^4)\right|_{x_k(t_0)}<0,$$
which is valid on $U\times [t_0, t_0+\de)$ as well, where $U$ is an open neighborhood of $x_k(t_0)$
in $M_{t_0}$. For $t_1$ is sufficiently close to $t_0,$ there exists $x_k(t_1)\in U\times [t_0, t_0+\de),$
such that $x_k(t_1)$ is evolved from some point $y$ on $M_{t_0}.$ Therefore,
$$||B||^2(x_k(t_1))<||B||^2(y)<\sup_{M_{t_0}}||B||^2.$$
Let $k\to\infty,$ we have the desired estimate in $[t_0, t_0+\de'),$ hence in $[0, T).$
\end{proof}

\bigskip

\subsection{Gauss Maps under the Evolution}

For any $p\in M$ let $\{e_1,\cdots, e_m\}$ be a local orthonormal frame field near $p$.
Define the Gauss map $\g: p\to \g(p)$ which is obtained by parallel translation of
$T_pM$ to the origin in the ambient space $\ir{m+n}_n.$ The image of the Gauss map
lies in a pseudo-Grassmannian $\grs{m}{n}^n$- the totality of all the space-like $m$-planes
in $\ir{m+n}_n$. It is a specific Cartan-Hadamard manifold.

For any $P\in \grs{m}{n}^n$, there are $m$ vectors $v_1, \cdots, v_m$ spanning $P$. Then
we have Pl\"ucker coordinates $v_m\wedge\cdots\wedge v_m$ for $P$ up to a constants.
The Gauss map $\g$ can be described by $p\to e_1\wedge\cdots\wedge e_m.$ Since
\begin{equation*}\begin{split}
d(e_1\wedge\cdots\wedge e_m)&=d e_1\wedge\cdots\wedge e_m+\cdots + e_1\wedge\cdots\wedge d e_m\\
&=\om_{\a 1}e_\a\wedge e_2\wedge\cdots\wedge e_m +\cdots+e_1\wedge\cdots\wedge e_{m-1}\wedge\om_{\a m}e_\a\\
&=\om_{\a i}e_{\a i}
\end{split}\end{equation*}
and the canonical metric on $\grs{m}{n}^n$ is defined by
$$ds^2 = \sum_{\a, i}\om_{\a i}^2,$$
where $e_{\a i}=e_1\wedge\cdots\wedge e_{i-1}\wedge e_\a\wedge e_{i+1}\wedge\cdots\wedge e_m$
are orthonomal basis for $T\grs{m}{n}^n.$ It follows that
$$\g^*\om_{\a i} = h_{\a ij}\om_j$$
which means that the relation
$$e(\g)=\frac 12 ||B||^2.$$
The tension field of the Gauss map
\begin{equation}\begin{split}
\tau(\g)&=h_{\a ijj}e_{\a i}=h_{\a jji}e_{\a i}\\
&=h_{\a jji}e_1\wedge\cdots\wedge e_{i-1}\wedge e_\a\wedge e_{i+1}\wedge\cdots\wedge e_m\\
&=\sum_i e_1\wedge\cdots\wedge e_{i-1}\wedge \n_{e_i}H\wedge e_{i+1}\wedge\cdots\wedge e_m,\label{TG}
\end{split}\end{equation}
where we use the Codazzi equation.

Now, we compute the evolution of the Gauss map under the mean curvature flow. From
$$\g(t)=\frac 1{\sqrt{g}}\,e_1(t)\wedge\cdots \wedge e_m(t)$$
and (\ref{tv}), we have
\begin{equation}\begin{split}
\dt{\g}&=-\frac 1g\dt{\sqrt{g}} \left(e_1(t)\wedge\cdots\wedge e_m(t)\right)\\
&\qquad +\frac 1{\sqrt{g}}\left(\dt{e_1(t)}\wedge\cdots\wedge e_m(t)+\cdots e_1(t)\wedge\cdots\wedge\dt{e_m(t)}\right)\\
&=\frac{|H|^2}{\sqrt{g}}\left(e_1(t)\wedge\cdots\wedge e_m(t)\right)\\
&\qquad+\frac 1{\sqrt{g}}\left(\bn_{e_1}H\wedge\cdots\wedge e_m(t)+\cdots+e_1(t)\wedge\cdots\wedge\bn_{e_m}H\right)\\
&=\frac{|H|^2}{\sqrt{g}}\left(e_1(t)\wedge\cdots\wedge e_m(t)\right)\\
&\qquad+\frac 1{\sqrt{g}}\left<\bn_{e_1}H, e_j(t)\right>e_k(t) g^{jk}\wedge e_2(t)\wedge\cdots\wedge e_m(t)+\cdots\\
&\qquad+\frac 1{\sqrt{g}}e_1(t)\wedge\cdots\wedge\left<\bn_{e_m}H, e_j(t)\right>e_k(t)g^{jk}\\
&\qquad+\frac 1{\sqrt{g}}\left(\n_{e_1}H\wedge\cdots\wedge e_m(t)+\cdots+e_1(t)\wedge\cdots\wedge\n_{e_m}H\right).
\end{split}\end{equation}
Since
\begin{equation}\begin{split}
&\left<\bn_{e_1}H, e_j(t)\right>e_k(t) g^{jk}\wedge e_2(t)\wedge\cdots\wedge e_m(t)+\cdots
e_1(t)\wedge\cdots\wedge\left<\bn_{e_m}H, e_j(t)\right>e_k(t)g^{jk}\\
&=-\left(\left<H,\bn_{e_1}e_j(t)\right>g^{1j}+\left<H,\bn_{e_m}e_j(t)\right>g^{jm}\right)
e_1(t)\wedge\cdots\wedge e_m(t)\\
&\qquad =-|H|^2\left(e_1(t)\wedge\cdots\wedge e_m(t)\right),
\end{split}\end{equation}
we have
\begin{equation}
\dt{\g}=\frac 1{\sqrt{g}}
\left(\n_{e_1}H\wedge\cdots\wedge e_m(t)+\cdots+e_1(t)\wedge\cdots\wedge\n_{e_m}H\right).\label{tG}
\end{equation}
We may assume $\{e_1(t_0),\cdots,e_m(t_0)\}$ form an orthonormal basis of $F(M)$ at $(p, t_0)$.
Then, from (\ref{TG}) and (\ref{tG}) we obtain the following equation which is the Lorentzian version
of a result in \cite{wang2}.
\begin{thm}\label{hfg}
\begin{equation}
\dt{\g}=\tau(\g(t)).\label{prv}
\end{equation}
\end{thm}

\bigskip
\subsection{Proof of the First Main Theorem}

We study the Gauss image under the flow.
The relevant Bernstein type theorem inspires us to consider the bounded Gauss image of the initial submanifold.
In fact, any geodesic ball in any Cartan-Hadamard manifold is convex.
Precisely, we have the following "confinable property" of the Gauss image under the mean curvature flow.
\begin{thm} \label{pbg}
Let $M$ be a complete space-like  $m-$submanifold in $\ir{m+n}_n$ with bounded curvature.
If the image under Gauss map is contained  in a bounded geodesic ball in $\grs{m}{n}^n,$ then the images of all
the submanifolds under the MCF are also contained in the same geodesic ball.
\end{thm}
\begin{proof}
Let $h$ be any function on $\grs{m}{n}^n.$ The composition function $h\circ \g$ of $h$
with the Gauss map $\g$ defines a function on $M_t=F(M, t)$. We have
$$\dt{}(h\circ\g)=dh\left(\dt{\g}\right)=dh(\tau(\g)).$$
By the composition formula (see \cite{X2}, p.28)
$$\De(h\circ\g)=\text{Hess}(h)(\g_*e_i,\g_*e_i)+dh(\tau(\g)),$$
where $\{e_i\}$ is a local orthonormal frame field on $M_t$. It follows that
\begin{equation}
\left(\dt{}-\De\right)h\circ\g=-\text{Hess}(h)(\g_*e_i,\g_*e_i).
\end{equation}
Noting $\grs{m}{n}^n$ has non-positive sectional curvature, the standard Hessian comparison theorem
implies
$$\text{Hess}(\tilde r)\ge\frac 1{\tilde r}(\tilde g-d\tilde r\otimes d\tilde r),$$
where $\tilde r$ is the distance function from a fixed point in $\grs{m}{n}^n$,
$\tilde g$ is the metric tensor on $\grs{m}{n}^n$. Choose $h=\tilde r^2$ and we have
$$\text{Hess}(h)\ge 2\tilde g.$$
Hence,
$$\text{Hess}(h)(\g_*e_i, \g_*e_i)\ge 4\,e(\g)=2\,||B||^2.$$

On the other hand,
\begin{equation}\begin{split}
|\n(h\circ\g)|^2&=\left<\n h, \g_*e_i\right>\left<\n h, \g_*e_i\right>\\
&=\left<2\tilde r\n \tilde r,\g_*e_i\right>\left<2\tilde r\n \tilde r,\g_*e_i\right>\\
&\le 8\tilde r^2\,e(\g)=4\tilde r^2||B||^2.
\end{split}\end{equation}
Therefore
\begin{equation}
\left(\dt{}-\De\right)h\circ\g\le - 2\,||B||^2\le -\frac 1{2\,h\circ \g}|\n(h\circ\g)|^2.
\label{EGD}\end{equation}

First of all, by Theorem \ref{cur-es} we always have bounded curvature for the smooth solution of
(\ref{mcf}). Denote $u=h\circ \g.$ Then, by (\ref{EGD})
$$\left(\dt{}-\De\right)u\le 0.$$
Let $u_k(t_0)$ be a sequence of points on $M_{t_0}$ such that when $k\to \infty$
$$u(x_k(t_0))\to\sup_{M_0} u>0$$
and when $k$ is sufficiently large
\begin{equation*}\begin{split}
|\n u|(x_k(t_0))<\ep,\\
\De u (x_k(t_0))<\ep
\end{split}
\end{equation*}
for any $\ep>0$. Define
$$u_1=(u-\sup_{M_{t_0}}u)-\de(t-t_0)-\de$$
for any $\de>0$ and when $t=t_0$
$$u_1\le -\de<0.$$
Thus,
$$\left.\dt{u_1}\right|_{x_k(t_0)}=\left.\dt{u}\right|_{x_k(t_0)}-\de\le\left.\De u\right|_{x_k(t_0)}-\de<0.$$
The the similar argument as that in the proof of Theorem \ref{cur-es} gives
$$u_1(x_k(t_1))<u_1(y)\le -\de$$
for $t_1$ close to $t_0$, namely
$$u(x_k(t_1))\le\sup_{M_{t_0}}u.$$
Let $k\to\infty$ gives
$$\sup_{M_{t_1}}u\le\sup_{M_{t_0}}u.$$
\end{proof}

Choose a Lorentzian base $\{\ep_i, \ep_\a\}$ in $\ir{m+n}_n$ with space-like
$\{\ep_i\}$ and time-like $\{\ep_\a\}.$ For a space-like submanifold
$F:M\to\ir{m+n}_n$ we assume $0\in M$ and define coordinate functions
$$x^i=\left<F, \ep_i\right>,\;y^\a=\left<F, \ep_\a\right>.$$
Denote
$$x=\sqrt{\sum_{i=1}^m (x^i)^2},\quad  y=\sqrt{\sum_{\a=m+1}^{m+n} (y^\a)^2},$$
such that $|F|^2=x^2-y^2.$ It is non-negative. The function $y$ is called the height function of $M$.
The growth
$$y^2\le x^2$$
is always satisfied for any space-like submanifold.
\begin{thm}
Let $F:M\to \ir{m+n}_n$ be a space-like complete $m-$submanifold which has bounded curvature
and bounded Gauss image. Then the mean curvature flow equation (\ref{mcf}) has long time smooth solution.
\end{thm}
\begin{proof}
Let $p:\ir{m+n}_n\to \ir{m}$ be the natural projection defined by
$$p(x^1,\cdots,x^m;x^{m+1},\cdots,x^{m+n})=(x^1,\cdots,x^m),$$
which induces a map from $M$ into $\ir{m}.$ It is a smooth map from a complete manifold to $\ir{m}.$
For any vector $v=(v^1,\cdots, v^m,v^{m+1},\cdots, v^{m+n})$ tangent to $M\subset \ir{m+n}_n$,
we define
$$\left<p_*v, p_*v\right>_{\ir{m}}=\left<v, v\right>+\sum_\a(v^{m+\a})^2\ge \left<v, v\right>.$$
This means that $p$ increases the distance. It follows that $p$ is a covering map, and
a deffeomorphism, since $\ir{m}$ is simply connected. Hence, the induced Riemannian metric on
$M$ can be expressed as $(\ir{m}, ds^2)$ with
$$ ds^2=g_{ij}dx^i\,dx^j.$$
In this chart, the identity map $(\ir{m}, ds^2)\to (\ir{m}, ds_0^2)$ is a distance
increasing map, where $ds_0^2$ is the Euclidean metric on $\ir{m}$. It follows that any eigenvalue
of $(g_{ij})$ is not big than $1$.

Choose $P_0$ as an $m-$plane spanned by
$\ep_1\wedge\cdots\wedge \ep_m.$ At each point in $M$ its image
$m-$plane $P$ under the Gauss map is spanned by
$$f_i=\ep_i+ c_{is}\ep_{m+s}$$
and
$$\sqrt{g}=|f_1\wedge\cdots\wedge f_m|.$$
We then have
\begin{equation}\begin{split}
\left<P,P_0\right>&=\frac 1{\sqrt{g}}\text{det}(\left<f_i,\ep_j\right>)\cr
 &=g^{-\frac 12}.
\end{split}\label{e-p}\end{equation}

Now, drawing a minimal geodesic $C(r)$ between $P_0$ and $P$ parameterized
by arc length $r$. By a result in  \cite{Wong}, $C(r)$ can be represented
by $P(r)$ which is spanned by
$$h_i= \ep_i + z_{is}(r)\ep_{m+s},$$
where
$$z_{is}(r)=
\begin{pmatrix} \tanh(\la_1 r)&&0\cr
          &\ddots&&0\cr
      0&&\tanh(\la_m r)&&
  \end{pmatrix}$$
for $\sum_i\la_i^2=1.$ Let

$$\tilde h_1=\cosh(\la_1 r)h_1,\; \cdots,\; \tilde h_m=\cosh(\la_m r)h_m.$$
Since
\begin{equation*}\begin{split}
|h_i|^2
&=\left<\ep_i+z_{is}(r)\ep_{m+s}, \ep_i+z_{it}(r)\ep_{m+t}\right>\cr
 &= 1- \tanh^2(\la_i r) = \frac 1{\cosh^2(\la_i r)},
\end{split}\end{equation*}
the vectors $\tilde h_1,\cdots,\tilde h_m$ are orthonormal.
Therefore, we can compute the inner product $\left<P_0,P\right>$
again by
$$
\left<P_0,P\right>=\text{det}\left(\left<\ep_i,\tilde
h_j\right>\right) =\prod_{i=1}^m\cosh(\la_i r).$$
If the distance $r$ between $P_0$ and $P$ is bounded by a
finite number $R$, then combine the above formula and (\ref{e-p}) yields
$$\sqrt{g}\ge (\prod_{i=1}^m\cosh(\la_i R))^{-1}.$$
Thus, we prove that any eigenvalue of the induced metric tensor of a complete
space-like $m-$submanifold in $\ir{m+n}_n$ with bounded Gauss image is uniformly bounded.
Noting Theorem \ref{pbg}, we know that the equation (\ref{mcf}) is uniformly parabolic
and has a unique smooth solution on some short time interval. By the curvature estimate
(see Theorem \ref{cur-es}), we have uniform estimate on $||B||$.
Then we can proceed as  in \cite{h2} (Prop. 2.3) to estimate all derivatives of
$B$ in terms of their initial data
$$\sup_{M_t}||\n^qB||\le C(m),$$
where $C(m)$ only depends on $q, m$ and $\sup_{M_0}||\n^jB||$ for
$0\le j\le q.$ It follows that this solution can be extended to 
all $t>0.$
\end{proof}

It is easy to verified that
$$\left(\dt{}-\De\right) y^\a=0$$
and
$$\left(\dt{}-\De\right) y^2= -2\sum|\n y^\a|^2\le 0.$$
Omori-Yau maximum principle implies that if $M_0$ has finite curvature and finite height
function, then the height function of $M_t$ is also finite under the evolution.

By (\ref{ENB}) and (\ref{EGD}) we have
$$\left(\dt{}-\De\right)(2 t ||B||^2+ h\circ\g)\le -\frac{4t}{n}||B||^4,$$
where $h$ denotes the square of the distance function on $\grs{m}{n}^n$, the target
manifold of the Gauss map. By using  maximum  principle again we have an estimate for $t>0$
$$\sup_{M_t}||B||^2\le \frac ct,$$
where $c$ is a constant depending on the bound of the Gauss image of the initial submanifold.

If the height function is going to infinity, we can consider rescaled mean curvature flow
as done by Ecker-Huisken \cite{E-H1}. Define
$$\tilde F(\tilde t)=\frac 1{\sqrt{2t+1}}F(t),$$
where
$$\tilde t=\log(2t+1).$$
Hence
$$\frac {\p}{\p \tilde t}\tilde F=\tilde H-\tilde F.$$
It is not hard to verify that the Gauss map $\tilde \g$ of the rescaled mean curvature flow is
as same as the original $\g.$ Furthermore, the previous estimates translate to
$$|\tilde A|^2\le (2 t+1)|A|^2\le C$$
which is dependent on the initial bound on $M.$

Choose a Lorentzian frame field near $p\in M$ along $M$ in $\ir{m+n}_n$, such that
$e_i\in TM$ and $e_\a\in NM$ with $\n_{e_j}e_i|_p=\n_{e_i}e_\a|_p=0.$ We have

\begin{lem}
\begin{equation*}
\left(\dt{}-\De\right)\left<F, e_\a\right>=2\,\left<H, e_\a\right>-S_{\a\be}\left<F, e_\be\right>
+C_{\a\be}\left<F, e_\be\right>
\end{equation*}
with anti-symmetric $C_{\a\be}$  in $\a$ and $\be$, and
\begin{equation}\begin{split}
\left(\dt{}-\De\right)\sum_\a\left<F, e_\a\right>^2
&=4\,\left<H,e_\a\right>\left<F, e_\a\right>-2\,\sum_\a|\n\left<F, e_\a\right>|^2
-2\,S_{\a\be}\left<F, e_\a\right>\left<F, e_\be\right>\\
&\le C\left(\sum_\a\left<F, e_\a\right>^2+1\right)-2\,\sum_\a|\n\left<F, e_\a\right>|^2.
\end{split}\end{equation}
\end{lem}
\begin{proof}
Since at the point $p$
\begin{equation}
\begin{split}
\n_{e_i}A^{e_\a}(e_i)&=\n_{e_i}\left<B_{ij},e_\a\right>e_j\\
&=\n_{e_i}\left(\left<\bn_{e_j}e_i,e_\a\right>e_j\right)\\
&=\left<\bn_{e_j}\bn_{e_i}e_i, e_\a\right>e_j\\
&=\left<\bn_{e_j}(\n_{e_i}e_i+B_{ii}), e_\a\right>e_j\\
&=\left<B_{j\n_{e_i}e_i}, e_\a\right>+\left<\bn_{e_j}H, e_\a\right>e_j\\
&=\left<\n_{e_j}H, e_\a\right>e_j,
\end{split}
\end{equation}
\begin{equation}
\begin{split}
\De\left<F, e_\a\right>&=e_i e_i\left<F, e_\a\right>\\
&=-e_i\left<F, A^{e_\a}(e_i)\right>\\
&=-\left<e_i, A^{e_\a}(e_i)\right>-\left<F,\bn_{e_i}A^{e_\a}(e_i)\right>\\
&=-\left<H, e_\a\right>-\left<F, \n_{e_i}A^{e_\a}(e_i)\right>-\left<F, B_{e_iA^{e_\a}(e_i)}\right>\\
&=-\left<H, e_\a\right>-\left<\n_{e_i}H, e_\a\right>\left<F, e_i\right>-\left<F, B_{e_iA^{e_\a}}(e_i)\right>\\
&=-\left<H, e_\a\right>-\left<\n_{e_i}H, e_\a\right>\left<F, e_i\right>-\left<F,B_{ij}\right>\left<B_{ij},e_\a\right>\\
&=-\left<H, e_\a\right>-\left<\n_{e_i}H, e_\a\right>\left<F, e_i\right>+S_{\a\be}\left<F, e_\be\right>.
\end{split}
\end{equation}
On the other hand,
$$\left<\dt{e_\a}, e_i\right>=-\left<e_\a, \dt{e_i}\right>=-\left<e_\a,\n_{e_i}H\right>,$$
$$\dt{e_\a}=-\left<e_\a,\n_{e_i}H\right>e_i+C_{\a\be}e_\be$$
with anti-symmetric $C_{\a\be}$  in $\a,\, \be.$
It follows that
\begin{equation}\begin{split}
\dt{}\left<F, e_\a\right>&=\left<H, e_\a\right>+\left<F, \dt{e_\a}\right>\\
&=\left<H, e_\a\right>-\left<e_\a, \n_{e_i}H\right>\left<F, e_i\right>+C_{\a\be}\left<F, e_\be\right>.
\end{split}
\end{equation}
Furthermore, we have
\begin{equation}\begin{split}
\dt{}\sum_\a\left<F, e_\a\right>^2&=2\,\sum_\a\left<F, e_\a\right>\dt{}\left<F, e_\a\right>\\
&=2\,\sum_\a\left<H, e_\a\right>\left<F, e_\a\right>
-2\,\sum_\a\left<\n_{e_i}H, e_\a\right>\left<F, e_i\right>\left<F, e_\a\right>\\
&\hskip2in +2\,C_{\a\be}\left<F, e_\be\right>\left<F, e_\a\right>\\
&=2\,\sum_\a\left<H, e_\a\right>\left<F, e_\a\right>
-2\,\sum_\a\left<\n_{e_i}H, e_\a\right>\left<F, e_i\right>\left<F, e_\a\right>,
\end{split}
\end{equation}
and
\begin{equation}\begin{split}
\De\left<F, e_\a\right>^2&=2\,|\n \left<F, e_\a\right>|^2+2\,\left<F, e_\a\right>\De\left<F, e_\a\right>\\
&=2\,|\n \left<F, e_\a\right>|^2
+2\,\left<F, e_\a\right>(-\left<H, e_\a\right>-\left<\n_{e_i}H, e_\a\right>\left<F, e_j\right>
+S_{\a\be}\left<F, e_\be\right>).
\end{split}
\end{equation}
Hence,
\begin{equation}\begin{split}
\left(\dt{}-\De\right)\sum_\a\left<F, e_\a\right>^2&=4\,\left<H, e_\a\right>\left<F, e_\a\right>
-2\,\sum_\a|\n \left<F, e_\a\right>|^2-2\,S_{\a\be}\left<F, e_\a\right>\left<F, e_\be\right>.
\end{split}
\end{equation}
Noting the estimates of $||H||^2$ , $||B||^2$ and $S_{\a\be},$ we obtain the desired estimate.
\end{proof}
The subsequent estimates in \cite {E-H1} can be carried out in the same way to derive
\begin{thm}
Suppose that $F:M\to \ir{m+n}_n$ is a space-like complete $m$-submanifold with bounded curvature and
bounded Gauss image. If in addition assume that
$$\sum_\a\left<F, e_\a\right>^2\le C'(1+|F|^2)^{1-\de}$$
is valid on $M$ for some constants $C'<\infty, \de>0,$ then the solution $\tilde M_{\tilde t}$ of the
rescaled equation converges for $\tilde t \to \infty$ to a limiting submanifold $\tilde M_\infty$
satisfying the equation
$$F^\perp =H.$$
\end{thm}

\Section{Submanifolds in Euclidean space}{Submanifolds in Euclidean space}

Let $F:M\to \ir{m+n}$ be an $m-$submanifold in $(m+n)-$dimensional Euclidean space
with the second fundamental form $B$ which can be viewed as a cross-section of
the vector bundle Hom($\odot^2TM, NM$) over $M,$ where  $TM$ and
$NM$ denote the tangent bundle and the normal bundle  along $M$, respectively.
A connection on Hom($\odot^2TM, NM$) can be induced from those of $TM$
and $NM$ naturally. We investigate the higher codimension $n\ge 2$ situation in this section.

For $\nu\in\G(NM)$ the shape operator $A^\nu: TM\to TM$  satisfies
$$\left<B_{X Y}, \nu\right>=\left<A^\nu(X), Y\right>.$$

The second fundamental form, curvature tensors of the submanifold, curvature tensor
of the normal bundle and that of the ambient manifold satisfy the Gauss equations,
the Codazzi equations and the Ricci equations.

Taking the trace of $B$ gives the mean curvature vector
$H$ of $M$ in $\ir{m+n}$,  a cross-section of the normal bundle.

Choose a local orthonormal frame field $\{e_i,e_\a\}$ along $M$ with dual frame field
$\{\omega_i,\omega_\a\}$, such that $e_i$ are tangent vectors to $M$.
The induced Riemannian metric of $M$ is given by $ds_M^2 =\sum_i\omega_i^2$
and the induced structure equations of $M$ are
\begin{equation*}\begin{split}
   & d\omega_i = \omega_{ij}\wedge\omega_j,\qquad
                 \omega_{ij}+\omega_{ji} = 0,\cr
   &d\omega_{ij}= \omega_{ik}\wedge\omega_{kj}+\omega_{i\a}\wedge\omega_{\a j},\cr
   &\Omega_{ij} = d\omega_{ij}-\omega_{ik}\wedge\omega_{kj}
                = -\frac 12 R_{ijkl}\omega_k\wedge\omega_l.
\end{split}\end{equation*}
By Cartan's lemma we have
$$\omega_{\a i} = h_{\a ij}\omega_j.$$

\bigskip

\subsection{Evolution Equation under the Mean Curvature Flow}

As the same derivation as in the space-like $m-$submanifold in pseudo-Euclidean space, we have the following formula.

\begin{pro}
\begin{equation}\begin{split}
(\n^2 B)_{XY}
&=\n_X\n_Y H+ \left<B_{Xe_i},H\right>B_{Ye_i}-\left<B_{XY},B_{e_ie_j}\right>B_{e_ie_j}\\
&+2\left<B_{Xe_j},B_{Ye_i}\right>B_{e_ie_j}-\left<B_{Ye_i},B_{e_ie_j}\right>B_{Xe_j}
-\left<B_{Xe_i},B_{e_ie_j}\right>B_{Ye_j}.\label{LB}
\end{split}\end{equation}
\end{pro}

Denote

$$B_{ij}=B_{e_ie_j}=(\bar\n_{e_i}e_j)^N=h_{\a ij}e_\a,$$
where $\{e_\a\}$ is a local orthonormal frame field of the normal bundle near
$x\in M.$ Let $S_{\a\be}=h_{\a ij}h_{\be ij}$. Then $|B|^2=\sum_{\a}S_{\a\a}$.

Noting
\begin{equation*}
-\left<B_{kl},B_{ij}\right>\left<B_{ij},B_{kl}\right>
=-h_{\a kl}h_{\a ij}h_{\be ij}h_{\be kl}=-\sum_{\a,\be} S_{\a\be}^2,
\end{equation*}
\begin{equation*}\begin{split}
&2\,\left<B_{il},B_{jk}\right>\left<B_{kl},B_{ij}\right>
  -2\,\left<B_{jk},B_{kl}\right>\left<B_{il},B_{ij}\right>\\
&=2\sum_{\a\ne\be}
\left(\left<A^{e_\be}A^{e_\a}, A^{e_\a}A^{e_\be}\right>-2\left<A^{e_\be}A^{e_\a}, A^{e_\be}A^{e_\a}\right>\right)\\
&=-\sum_{\a\ne\be}|[A^{e_\a}, A^{e_\be}]|^2,
\end{split}\end{equation*}
we then have
$$\left<\n^2B,B\right>=\left<\n_i\n_jH,B_{ij}\right>
                       +\left<B_{ik},H\right>\left<B_{il},B_{kl}\right>
                -\sum_{\a\ne\be}|[A^{e_\a}, A^{e_\be}]|^2-\sum_{\a,\be}S_{\a\be}^2.$$
The following expression follows immediately.
\begin{pro}
\begin{equation}\begin{split}
\De|B|^2=2\,|\n B|^2+2\,\left<\n_i\n_jH,B_{ij}\right>
       &+\,2\left<B_{ij},H\right>\left<B_{ik},B_{jk}\right>\\
       &\quad -\,2\,\sum_{\a\ne\be}|[A^{e_\a}, A^{e_\be}]|^2-2\,\sum_{\a,\be}S_{\a\be}^2\label{NLB'}
\end{split}\end{equation}
\end{pro}

We now consider the MCF for a submanifold  in $\ir{m+n}.$ Namely,
consider a one-parameter family $F_t=F(\cdot, t)$ of immersions $F_t:M\to \ir{m+n}$
with corresponding images $M_t=F_t(M)$ such that
\begin{equation}\begin{split}
\dt{}F(x, t)&=H(x, t),\qquad x\in M\\
F(x, 0)&=F(x)
\end{split}\label{mcf1}
\end{equation}
is satisfied, where $H(x, t)$ is the mean curvature vector of $M_t$ at $F(x, t)$ in $\ir{m+n}$
We also have
\begin{equation}
\dt{g_{ij}}=-2\left<H,B_{ij}\right>.\label{tg1},
\end{equation}
\begin{equation}
\dt{g^{ij}}=2\,g^{ik}g^{jl}\left<H,B_{kl}\right>\label{tg1'}
\end{equation}
and
\begin{equation}
\dt{g}=-2\,|H|^2\,g.\label{tv1},
\end{equation}
where $g=\det(g_{ij}).$ We now derive the evolution equation for the squared norm of the second fundamental
form.

\begin{lem}
The second fundamental form  satisfies
\begin{equation}
\left(\dt{}-\De\right)|B|^2\le - \,2\,|\n|B||^2 + 3|B|^4\label{ENB1},
\end{equation}
\end{lem}
\begin{rem}
Compare (\ref{ENB1}) with (\ref{ENB}), we see that now the curvature estimates is not
directly as that in the last section.
\end{rem}
\begin{proof}
For fixed $x_0,\,t_0$ choose a local orthonormal frame $\{e_i\}$ of $M_{t_0}$ near $x_0$
which is normal at $x_0$. By the immersion $F_{t_0}$ we have $\{e_i\}$ on $M$, which is not orthonormal in general.
Then by $F_t$ we obtain $\{F_{t*}e_i\}$ which is denoted by $\{e_i\}$ for simplicity.
We also choose a local orthonormal frame field $\{e_\a\}$ of the normal bundle of $M_t$ near $x_0$.
Then at $(x_0, t_0)$
\begin{equation}\begin{split}
\dt{h_{\a ij}}&= \bn_{\dt{}}\left<\bn_{e_i}e_j,e_\a\right>\\
&=\left<\bn_H \bn_{e_i}e_j,e_\a\right>+\left<\bn_{e_i}e_j,\bn_He_\a\right>\\
&=\left<\bn_{e_i}\bn_{e_j}H,e_\a\right>+\left<B_{ij},\bn_He_\a\right>\\
&=\left<\bn_{e_i}\left(\n_{e_i}H+\left(\bn_{e_j}H\right)^T\right),e_\a\right>
+\left<B_{ij},\bn_He_\a\right>\\
&=\left<\n_{e_i}\n_{e_j}H,e_\a\right>-h_{\a ik}h_{\be jk}H_\be
     +h_{\be ij}\left<\bn_He_\a,e_\be\right>.\label{TB1'}
\end{split}\end{equation}

Since in a non-orthonormal frame field $g_{ij}=\left<F_*e_i. F_*e_j\right>$ (except at $t_0$)
is not a unit matrix,
$$|B|^2= g^{ik}g^{jl}h_{\a i j}\,h_{\a k l}.$$
We have at $(x_0, t_0)$
\begin{equation}
\dt{|B|^2}=2\,\dt{g^{ik}}h_{\a i j}h_{\a k j}+2\,\dt{h_{\a i j}}h_{\a i j}.
\label{TB2'}
\end{equation}
From (\ref{TB1'}) we have
\begin{equation}
\dt{h_{\a i j}}h_{\a i j}
=h_{\a i j}\left<\n_{e_i}\n_{e_j}H, e_\a\right>
  -h_{\a ij}h_{\a ik}h_{\be jk}H_\be
\label{TB3'}
\end{equation}
Noting (\ref{tg1'}), we have
\begin{equation}
\dt{g^{ik}}h_{\a i j}h_{\a k j}
=2\,h_{\a ij}h_{\a k j}\left<H,B_{ik}\right>
=2\,h_{\a ij}h_{\a k j}h_{\be ik}H_\be\label{TB4'}
\end{equation}
Substituting (\ref{TB3'}) and (\ref{TB4'}) into (\ref{TB2'})gives
\begin{equation}
\frac 12\dt{|B|^2}=\left<\n_i\n_jH,B_{ij}\right>
       +\left<B_{ij},H\right>\left<B_{ik},B_{jk}\right>\label{TB'}
\end{equation}
From (\ref{NLB'}) and (\ref{TB'}) we obtain
the evolution equation for the squared norm of the second fundamental form
\begin{equation}
\frac 12\left(\dt{}-\De\right)|B|^2=-|\n B|^2+\sum_{\a\ne\be}|[A^{e_\a}, A^{e_\be}]|^2+\sum_{\a, \be}S_{\a\be}^2.
\label{EB'}\end{equation}
We know from \cite{JSi} in general
$$\sum_{\a\ne\be}|[A^{e_\a}, A^{e_\be}]|^2+\sum_{\a, \be}S_{\a\be}^2\le \left(2-\frac 1n\right)|B|^4.$$
When the codimension $n\ge 2$ the above estimate was refined \cite{L-L}\cite{C-X}
$$\sum_{\a\ne\be}|[A^{e_\a}, A^{e_\be}]|^2+\sum_{\a, \be}S_{\a\be}^2\le\frac 32 |B|^4.$$
On the other hand, by  the Schwartz  inequality
$$\n |B|\le |\n B|.$$
Therefore, the inequality (\ref{ENB1}) is obtained.
\end{proof}

\subsection{Main Estimates}

For any $p\in M$ let $\{e_1,\cdots, e_m\}$ be a local orthonormal frame field near $p$.
Define the Gauss map $\g: p\to \g(p)$ which is obtained by parallel translation of
$T_pM$ to the origin in the ambient space $\ir{m+n}.$ The image of the Gauss map
lies in a Grassmannian $\grs{m}{n}$. It is a symmetric space of compact type.

For any $P\in \grs{m}{n}$, there are $m$ vectors $v_1, \cdots, v_m$ spanning $P$. Then
we have Pl\"ucker coordinates $v_m\wedge\cdots\wedge v_m$ for $P$ up to a constants.
The Gauss map $\g$ can be described by $p\to e_1\wedge\cdots\wedge e_m.$ Since
\begin{equation*}\begin{split}
d(e_1\wedge\cdots\wedge e_m)&=d e_1\wedge\cdots\wedge e_m+\cdots + e_1\wedge\cdots\wedge d e_m\\
&=\om_{\a 1}e_\a\wedge e_2\wedge\cdots\wedge e_m +\cdots+e_1\wedge\cdots\wedge e_{m-1}\wedge\om_{\a m}e_\a\\
&=\om_{\a i}e_{\a i}
\end{split}\end{equation*}
and the canonical metric on $\grs{m}{n}$ is defined by
$$ds^2 = \sum_{\a, i}\om_{\a i}^2,$$
where $e_{\a i}=e_1\wedge\cdots\wedge e_{i-1}\wedge e_\a\wedge e_{i+1}\wedge\cdots\wedge e_m$
are orthonomal basis for $T\grs{m}{n}$ (see \cite{X3}, pp. 188-194). It follows that
$$\g^*\om_{\a i} = h_{\a ij}\om_j$$
and the tension field of the Gauss map
\begin{equation}\begin{split}
\tau(\g)&=h_{\a ijj}e_{\a i}=h_{\a jji}e_{\a i}\\
&=h_{\a jji}e_1\wedge\cdots\wedge e_{i-1}\wedge e_\a\wedge e_{i+1}\wedge\cdots\wedge e_m\\
&=\sum_i e_1\wedge\cdots\wedge e_{i-1}\wedge \n_{e_i}H\wedge e_{i+1}\wedge\cdots\wedge e_m,\label{TG1}
\end{split}\end{equation}
where we use the Codazzi equation. In \cite{wang2}, there is the following relation.
\begin{pro}\label{wang2}
\begin{equation}
\dt{\g}=\tau(\g(t)).\label{prv}
\end{equation}
\end{pro}

We consider the mean curvature flow of a complete manifold.
We will assume that integration by parts is permitted and all integrals are finite
for the submanifolds and functions we will consider in the sequel.
We have the following {\sl maximum principle} for parabolic equations on complete manifolds.

Define the {\sl backward heat kernel} $\rho=\rho(x, t)$
by
$$\rho(x,t)=\frac{1}{(4\pi(t_0-t))^\frac{n}{2}}\,\exp\left({-\frac{|x|^2}{4(t_0-t)}}\right),\quad t_0>t,
\quad x\in\ir{m+n}.$$
We have the following formula. It is derived in the mean curvature flow in Euclidean space.
Since (\ref{tv1}), the formula is unchanged in higher codimension.
\begin{pro}(Huisken \cite{h2})
For a function $f(x,t)$ on $M$ we have
\begin{equation}\label{mon}
\dt{}\int_M f\rho d\mu_t=\int_M\left(\dt\,f-\Delta f\right)\rho d\mu_t
-\int_M f\rho\left|H+\frac{F^\perp}{2(t_0-t)}\right|^2d\mu_t,
\end{equation}
where $d\mu_t$ is the volume form of $M_t$.
\end{pro}

\begin{cor}(Ecker and Huisken \cite{E-H1})\label{max princ}
Suppose the function $f=f(x,t)$ satisfies the inequality
\begin{equation}\nonumber
\left(\dt{}-\Delta\right)f\le\langle {\bf a},\nabla f\rangle
\end{equation}
for some vector field $\bf a$ with uniformly bounded norm on $M\times[0,t_1]$ for some $t_1>0$, then
$$\sup_{M_t}f\le\sup_{M_0}f$$
for all $t\in[0,t_1]$.
\end{cor}

Now, we consider the convex  Gauss image situation which is preserved
under the flow, as shown in the following theorem.

\begin{thm}(\rm \ confinable property) \it
If the Gauss image of the initial submanifold $M$ is contained in a geodesic ball of
the radius $\rho_0<\frac{\sqrt{2}}{4}\pi$ in $\grs{m}{n}$, then the Gauss images of all the
submanifolds under the MCF are also contained in the same geodesic ball.
\label{cp}\end{thm}
\begin{proof}
We consider a smooth bounded function on $\grs{m}{n}$
$$h=1+\ep-\cos(\sqrt{2}\rho),$$
where $\rho$ is the distance function from a point in $\grs{m}{n},\,  \ep>0$ is a fixed constant.
When $\rho<\frac{\sqrt{2}}{4}\pi$, $h$ is convex. By the Hessian
comparison theorem we have
$$\text{Hess}(h)\ge 2\,\cos(\sqrt{2}\rho) \,g,$$
where $g$ is the metric tensor on $\grs{m}{n}.$ Hence,
$$\text{Hess}(h)(\g_*e_i, \g_*e_i)\ge 2\,\cos(\sqrt{2}\rho) \,|B|^2$$

The composition function $h\circ \g$ of $h$
with the Gauss map $\g$ defines a function on $M_t=F(M, t)$. We have
$$\dt{}(h\circ\g)=dh\left(\dt{\g}\right)=dh(\tau(\g)).$$
By the composition formula (see \cite{X2}, p.28)
$$\De(h\circ\g)=\text{Hess}(h)(\g_*e_i,\g_*e_i)+dh(\tau(\g)),$$
where $\{e_i\}$ is a local orthonormal frame field on $M_t$.

It follows that
\begin{equation}
\left(\dt{}-\De\right)h\circ\g\le -2\,\cos(\sqrt{2}\rho\circ\g)\,|B|^2.
\label{eg}\end{equation}
Thus, we can use Corollary  \ref{max princ} to get conclusion.
\end{proof}

For simplicity $h\circ \g$ is denoted by $h_1$ in the sequel. On 
the other hand,
\begin{equation}\begin{split}
|\n h_1|^2&=|\left<\n h, \g_*e_i\right>\left<\n h, \g_*e_i\right>|\\
&\le 2\,\sin^2(\sqrt{2}\rho\circ\g)|B|^2.
\end{split}\label{es}\end{equation}
From (\ref{eg}) and (\ref{es}) we have
\begin{equation}
\left(\dt{}-\De\right)h_1\le - \cos(\sqrt{2}\rho\circ\g) \,|B|^2
-\frac {\cos(\sqrt{2}\rho\circ\g)}{2\,\sin^2(\sqrt{2}\rho\circ\g)}|\n h_1|^2.
\label{EG}\end{equation}
For any $q>0$,
\begin{equation}\begin{split}
\left(\dt{}-\De\right)h_1^q
&=q\,h_1^{q-1}\left(\dt{}-\De\right)h_1-q\,(q-1)\,h_1^{q-2}|\n h_1|^2\\
&\le -q\, h_1^{q-1}\cos(\sqrt{2}\rho\circ\g)|B|^2\\
&\qquad-\left(q(q-1)h_1^{q-2}+q h_1^{q-1}\frac{\cos(\sqrt{2}\rho\circ\g)}{2\sin^2(\sqrt{2}\rho\circ\g)}\right)|\n h_1|^2.
\end{split}\label{QEG}\end{equation}
From (\ref{ENB1}) and(\ref{QEG}), we have
\begin{equation}\begin{split}
\left(\dt{}-\De\right)\left(|B|^2 h_1^q\right)
&=|B|^2\left(\dt{}-\De\right)h_1^q + h_1^q\left(\dt{}-\De\right)|B|^2-2\,\n|B|^2\cdot\n h_1^q\\
&\le (-q \cos(\sqrt{2}\rho\circ\g)+3h_1)|B|^4 h_1^{q-1}\\
&\qquad-\left[q(q-1)h_1^{q-2}+q\frac{\cos(\sqrt{2}\rho\circ\g)}{2\sin^2(\sqrt{2}\rho\circ\g)}h_1^{q-1}\right]
                                |B|^2|\n h_1|^2\\
&\qquad -2\,h_1^q|\n|B||^2-2\,\n |B|^2\cdot\n h_1^q\\
&=[3(1+\ep)-(3+q)\cos(\sqrt{2}\rho\circ\g)]h_1^{q-1}|B|^4\\
&\qquad -\left[q(q-1)+q\frac{\cos(\sqrt{2}\rho\circ\g)}{2\sin^2(\sqrt{2}\rho\circ\g)}h_1\right]h_1^{q-2}|B|^2|\n h_1|^2\\
 &\qquad -2\,h_1^q|\n|B||^2-2\,\n|B|^2\cdot\n h_1^q.
\end{split}\label{b2hq'}\end{equation}
By using the Young inequality we have
\begin{equation}\begin{split}
-2\n|B|^2\cdot\n h_1^q&=-(h_1^{-q}\n h_1^q)\cdot\n(|B|^2 h_1^q)+|B|^2 h_1^{-q}|\n h_1^q|^2\\
&\qquad -\n |B|^2\cdot\n h_1^q\\
&\le-q(h_1^{-1}\n h_1)\cdot\n(|B|^2 h_1^q)+q^2|B|^2 h_1^{q-2}|\n h_1|^2\\
&\qquad +\frac 12 q^2 h_1^{q-2}|B|^2|\n h_1|^2+2\,h_1^q|\n|B||^2\\
&\le-q(h_1^{-1}\n h_1)\cdot\n(|B|^2 h_1^q)+\\
&\qquad +\frac 32 q^2 h_1^{q-2}|B|^2|\n h_1|^2+2\,h_1^q|\n|B||^2.
\end{split}\label{b2hq1}\end{equation}
Thus, (\ref{b2hq'}) becomes
\begin{equation}\begin{split}
\left(\dt{}-\De\right)\left(|B|^2 h_1^q\right)
&\le [3(1+\ep)-(3+q)\cos(\sqrt{2}\rho\circ\g)]|B|^4 h_1^{q-1}\\
&\qquad + \left(\frac 12 q +1-\frac{\cos(\sqrt{2}\rho\circ\g)}{2\sin^2(\sqrt{2}\rho\circ\g)}h_1\right)
                                                                               q\, h_1^{q-2}|B|^2|\n h_1|^2\\
&\qquad\qquad -q(h_1^{-1}\n h_1)\cdot\n(|B|^2 h_1^q).
\end{split}\label{b2hq}\end{equation}

We now give the following result.
\begin{thm}
Let $M$ be a complete $m-$submanifold in $\ir{m+n}$ with bounded curvature.
Suppose that the image under the Gauss map from $M$ into $\grs{m}{n}$ lies in a geodesic ball of radius
$R_0<\frac{\sqrt{2}}{12}\pi.$ If $M_t$ is a smooth solution of (\ref{mcf1}), then there is
the following estimate
\begin{equation}
\sup_{M_t}|B|^2h_1^q\le \sup_{M_0}|B|^2h_1^q,
\label{ce}\end{equation}
where $q$ is a fixed constant depending on $R_0.$
\label{ces}\end{thm}
\begin{proof}
Let $r_0=\cos(\sqrt{2}R_0).$ Then $r_0>\frac{\sqrt{3}}{2}.$ It follows that
$$\frac {3}{2r_0}-\frac {r_0}{2(1-r_0^2)}<0.$$
It is possible to choose $\ep>0$ satisfying
\begin{equation}
\left(\frac {3}{2r_0}-\frac {r_0}{2(1-r_0^2)}\right)\ep+\frac{3}{2r_0}-\frac 12-\frac{r_0}{2(1+r_0)}\le 0.
\label{ep}\end{equation}
Set
$$q=3\left(\frac{1+\ep}{r_0}-1\right).$$
Then for $r=\cos(\sqrt{2}\rho\circ\g)\ge r_0$,
$$3(1+\ep)-(3+q)r=3(1+\ep)-3(1+\ep)\frac{r}{r_0}\le 0,$$
which implies the first term of the right hand side of (\ref{b2hq}) is non-positive.
Note
\begin{equation}\begin{split}
&\frac 12q+1-\frac{r}{2(1-r^2)}(1+\ep-r)\\
&=\frac 32\left(\frac{1+\ep}{r_0}-1\right)+1-\frac{r}{2(1-r^2)}(1+\ep-r)\\
&=\left(\frac {3}{2r_0}-\frac {r}{2(1-r^2)}\right)\ep+\frac{3}{2r_0}-\frac 12-\frac{r}{2(1+r)},
\end{split}\label{ep1}\end{equation}
which is non-increasing in $r$. Since (\ref{ep}), (\ref{ep1}) is non-positive when $r\ge r_0.$
It follows that under the conditions of the theorem, (\ref{b2hq}) becomes
\begin{equation}
\left(\dt{}-\De\right)\left(|B|^2 h_1^q\right)
\le  -q(h_1^{-1}\n h_1)\cdot\n(|B|^2 h_1^q).
\label{b2hq1}\end{equation}

From (\ref{es}) we have
$$|h_1^{-1}\n h_1|\le \frac{\sqrt{2}\sin(\sqrt{2}\rho\circ\g)}{1+\ep-\cos(\sqrt{2}\rho\circ\g)}|B|$$
Let
$$f(\th)=\frac{\sin\th}{1+\ep-\cos\th}.$$
Since $f''(\th)|_{f'(\th)=0}\le 0$,
\begin{equation}\begin{split}
f(\th)\le f(\th)|_{f'(\th)=0}&=\frac {\sqrt{1-\frac {1}{(1+\ep)^2}}}{1+\ep-\frac 1{1+\ep}}\\
&=\frac{\sqrt{(1+\ep)^2-1}}{(1+\ep)^2-1}=\frac{\sqrt{\ep(\ep+2)}}{\ep(\ep+2)}.
\end{split}\end{equation}
It follows that
$$|h_1^{-1}\n h_1|\le \frac{\sqrt{2\,\ep(\ep+2)}}{\ep(\ep+2)}|B|$$
Thus, we can use Corollary \ref{max princ} and the estimate (\ref{ce}) has been obtained.
\end{proof}

\begin{cor}
Suppose that the image under the Gauss map from $M$ into $\grs{m}{n}$ lies in a geodesic ball of radius
$R_0<\frac{\sqrt{2}}{12}\pi.$ If $M_t$ is a smooth solution of (\ref{mcf1}), then there is the following estimate
\begin{equation}
\sup_{M_t}|B|^2\le \frac ct,
\label{bde}\end{equation}
where $c$ is depends only on the bound of the Gauss image of its initial manifold.
\end{cor}
\begin{proof}
From (\ref{eg})
\begin{equation*}\begin{split}
\left(\dt{}-\De\right)h_1^q&=qh_1^{q-1}\left(\dt{}-\De\right)h_1-q(q-1)h_1^{q-2}|\n h_1|^2\\
&\le -2 qh_1^{q-1}\cos(\sqrt{2}\rho\circ\g)|B|^2-q(q-1)h_1^{q-2}|\n h_1|^2.
\end{split}\end{equation*}
It follows that
\begin{equation}\begin{split}
&\left(\dt{}-\De\right)(t|B|^2h_1^q+h_1^q)\le -q(h_1^{-1}\n h_1)\cdot\n(t|B|^2h_1^q)\\
&\qquad +|B|^2h_1^q-2 qh_1^{q-1}\cos(\sqrt{2}\rho\circ\g)|B|^2-q^2h_1^{q-2}|\n h_1|^2+qh_1^{q-2}|\n h_1|^2.
\end{split}\label{tb2hq}\end{equation}
Since
$$q(h_1^{-1}\n h_1)\cdot\n h_1^q=q^2h_1^{q-2}|\n h_1|^2,$$
(\ref{tb2hq}) becomes
\begin{equation*}\begin{split}
\left(\dt{}-\De\right)(t|B|^2h_1^q+h_1^q)&\le -q(h_1^{-1}\n h_1)\cdot\n(t|B|^2h_1^q+h_1^q)\\
&\qquad +|B|^2h_1^q-2 qh_1^{q-1}\cos(\sqrt{2}\rho\circ\g)|B|^2+qh_1^{q-2}|\n h_1|^2.
\end{split}\end{equation*}
Noting (\ref{es}), the above inequality becomes
\begin{equation}\begin{split}
\left(\dt{}-\De\right)(t|B|^2h_1^q+h_1^q)&\le -q(h_1^{-1}\n h_1)\cdot\n(t|B|^2h_1^q+h_1^q)\\
&\qquad +|B|^2h_1^{q-2}(h_1^2-2 qh_1\cos(\sqrt{2}\rho\circ\g)+2q\sin^2(\sqrt{2}\rho\circ\g)).
\end{split}\label{tb2hq1}\end{equation}

Let
\begin{equation}\begin{split}
A(r)&=(h_1^2-2 qh_1\cos(\sqrt{2}\rho\circ\g)+2q\sin^2(\sqrt{2}\rho\circ\g))\\
&=(1+\ep-r)^2 - 2q(1+\ep-r)r + 2q(1-r^2)=(1+\ep-r)^2-2q(r+\ep r-1),
\end{split}\end{equation}
where $r=\cos(\sqrt{2}\rho\circ\g)$. Since $A'(r)<0$ and $q=3\left(\frac {1+\ep}{r_0}-1\right)$,
then for $r\ge r_0$
$$A\le (1+\ep-r_0)^2-2q(r_0+\ep r_0-1)=(1+\ep-r_0)(\frac 6{r_0}-r_0-5-5\ep).$$
We know that $\ep$ is chosen by (\ref{ep}). If necessary we choose $\ep$ larger such that
$A\le 0.$ Therefore, from (\ref{tb2hq1}) we have
$$\left(\dt{}-\De\right)(t|B|^2h_1^q+h_1^q)\le -q(h_1^{-1}\n h_1)\cdot\n(t|B|^2h_1^q+h_1^q)$$
and by Corollary \ref{max princ} again we have the desired estimate
\end{proof}

\bigskip

\subsection{Proof of the Second Main Theorem}

We are now in a position to prove the following theorem.
\begin{thm}
Let $F:M\to \ir{m+n}$ be a  complete $m-$submanifold which has bounded curvature.
Suppose that the image under the Gauss map from $M$
into $\grs{m}{n}$ lies in a geodesic ball of radius $R_0<\frac{\sqrt{2}}{12}\pi.$
Then the mean curvature flow equation (\ref{mcf1}) has long time smooth solution.
\end{thm}
\begin{proof}
Let $P_0\in \grs{m}{n}$ be a fixed point which is described by
$$P_0=\ep_1\wedge\cdots\wedge\ep_m,$$
where $\ep_1,\cdots,\ep_m$ are orthonormal vectors in $\ir{m+n}$. Choose complementary
orthonormal vectors $\ep_{m+1},\cdots,\ep_{m+n}$, such that
$\{\ep_1,\cdots,\ep_m, \ep_{m+1},\cdots,\ep_{m+n}\}$
is an orhtonormal base in $\ir{m+n}$

Let $p:\ir{m+n}\to \ir{m}$ be the natural projection defined by
$$p(x^1,\cdots,x^m;x^{m+1},\cdots,x^{m+n})=(x^1,\cdots,x^m),$$
which induces a map from $M$ to $\ir{m}.$ It is a smooth map from a complete manifold to $\ir{m}.$

For any point $x\in M$ choose a local orhtonormal tangent frame field
$\{e_1,\cdots, e_m\}$ near $x$. Let $v=v_ie_i\in TM.$ Its projection
$$p_*v=\left<v_ie_i, \ep_j\right>\ep_j=v_i\left<e_i, \ep_j\right>\ep_j.$$

Now, we consider the case of the image under the Gauss map $\g$ containing in a
geodesic ball of radius $R_0<\frac {\sqrt{2}}{12} \pi$ and centered at $P_0$.
For any  $P\in\g(M)$,
$$w\mathop{=}\limits^{def.}\left<P, P_0\right>
=\left<e_1\wedge\cdots\wedge e_m, \ep_1\wedge\cdots\wedge \ep_m\right>=\det W,$$
where $W=(\left<e_i, \ep_j\right>).$ The Jordan angles between $P$ and $P_0$ are
$$\th_i=\cos^{-1}(\la_i),$$
where $\la_i^2$ are eigenvalues of the symmetric matrix $W^TW.$
It is well known that
$$W^TW=O^T\La O,$$
where $O$ is an orthogonal matrix and
$$\La=
\begin{pmatrix} \la_1^2&&0\cr
          &\ddots&\cr
      0&&\la_m^2
  \end{pmatrix},$$
where each $0<\la_i^2<1.$ We know that
$$w=\prod\cos\th_i.$$
On the other hand, the distance between $P_0$ and $P$ (see \cite{X3}, pp. 188-194)
$$d(P_0, P)=\sqrt{\sum\th_i^2}$$
which is less than $\frac{\sqrt{2}}{12}\pi$ by the assumption. It follows that
$$w>w_0=\left(\cos\frac{\sqrt{2}}{12}\pi\right)^m.$$
We now compare the length of any tangent vector $v$ to $M$ with its projection $p_*v.$
$$|p_*v|^2=\sum_{j=1}^m(v_i\left<e_i, \ep_j\right>)^2=(WV)^TWV,$$
where $V=(v^1, \cdots, v^m)^T$. Hence,
\begin{equation}
|p_*v|^2\ge (\la')^2|v|^2>w^2|v|^2>w_0^2|v|^2,
\label{di}\end{equation}
where $\la'=\min_i\{\la_i\}.$ The induced metric $ds^2$ on $M$ from $\ir{m+n}$ is complete,
so is the homothetic metric $\tilde ds^2=w_0^2ds^2.$ (\ref{di}) implies
$$p:(M, \tilde ds^2)\to(\ir{m},\text{canonical metric})$$
increases the distance. It follows that $p$ is a covering map from a complete manifold into $\ir{m}$ , and
a deffeomorphism, since $\ir{m}$ is simply connected. Hence, the induced Riemannian metric on
$M$ can be expressed as $(\ir{m}, ds^2)$ with
$$ ds^2=g_{ij}dx^i\,dx^j.$$
Furthermore, the immersion $F:M\to \ir{m+n}$ is realized by a graph $(x, f(x))$ with
$f:\ir{m}\to\ir{n}$ and
$$g_{ij}=\de_{ij}+\pd{f^\a}{x^i}\pd{f^\a}{x^j}$$
It follows that any eigenvalue of $(g_{ij})$ is not less than $1$.

At each point in $M$ its image $m$-plane $P$ under the Gauss map is spanned by
$$f_i=e_i+\pd{f^\a}{x^i}e_\a.$$
It follows that
$$|f_1\wedge\cdots\wedge f_m|^2=\text{det}\left(\de_{ij}
     +\sum_\a\pd{f^\a}{x^i}\pd{f^\a}{x^j}\right)$$
and
$$\sqrt{g}=|f_1\wedge\cdots\wedge f_m|.$$
The $m$-plane $P$ is also spanned by
$$p_i=g^{-\frac 1{2m}}f_i,$$
furthermore,
$$|p_1\wedge\cdots\wedge p_m|=1.$$
We then have
\begin{equation*}\begin{split}
\left<P,P_0\right>&=\text{det}(\left<\ep_i,p_j\right>)\cr
 &= \begin{pmatrix}
g^{-\frac 1{2m}}&&0\\
            &\ddots&\cr
    0&& g^{-\frac 1{2m}}\\
        \end{pmatrix}\\
  &=\frac 1{\sqrt{g}}>w_0
\end{split}\end{equation*}
and
$$\sqrt{g}\le \frac 1{w_0}$$

Thus, we prove that any eigenvalue of $(g_{ij})\le \frac 1{w_0^2}$.
Noting Theorem \ref{cp}, we know that the equation (\ref{mcf1}) is uniformly parabolic
and has a unique smooth solution on some short time interval. By the curvature estimate
(see Theorem \ref{cp} and Theorem \ref{ces}), we have uniform estimate on $|B|$.
Then we can proceed as  in \cite{h2} (Prop. 2.3) to estimate all derivatives of
$B$ in terms of their initial data
$$\sup_{M_t}|\n^qB|\le C(m),$$
where $C(m)$ only depends on $q, m$ and $\sup_{M_0}|\n^jB|$ for
$0\le j\le q.$ It follows that this solution can be extended to 
all $t>0.$
\end{proof}

We assume $0\in M$ and define coordinate functions
$$x^i=\left<F, \ep_i\right>,\;y^\a=\left<F, \ep_\a\right>.$$
Denote
$$x=\sqrt{\sum_{i=1}^m (x^i)^2},\quad  y=\sqrt{\sum_{\a=m+1}^{m+n} (y^\a)^2},$$
It is easy to verified that
$$\left(\dt{}-\De\right) y^\a=0$$
and
$$\left(\dt{}-\De\right) y^2= -2\sum|\n y^\a|^2\le 0.$$
Corollary \ref{max princ} implies that if the
height function of $M_0$ is finite, then the height function of $M_t$ is also
finite under the evolution.

If the height function is going to infinity, we can consider rescaled mean curvature flow
as done in \cite{E-H1}. Define
$$\tilde F(\tilde t)=\frac 1{\sqrt{2t+1}}F(t),$$
where
$$\tilde t=\log(2t+1).$$
Hence
$$\frac {\p}{\p \tilde t}\tilde F=\tilde H-\tilde F.$$
It is not hard to verify that the Gauss map $\tilde \g$ of the rescaled mean curvature flow is
as same as the original $\g.$ Furthermore, the previous estimates (\ref{bde}) translate to
$$|\tilde A|^2\le (2 t+1)|A|^2\le C$$
which is dependent on the initial bound on $M.$

We can  carried out in the same way as in \cite {E-H1} to derive

\begin{thm}
Let $F:M\to \ir{m+n}$ be a  complete $m$-submanifold with bounded curvature.
Suppose that the image under the Gauss map from $M$ into $\grs{m}{n}$ lies in a geodesic ball
of radius $R_0<\frac {\sqrt{2}}{12}\pi.$ If in addition assume that
$$\sum_\a\left<F, e_\a\right>^2\le C'(1+|F|^2)^{1-\de}$$
is valid on $M$ for some constants $C'<\infty, \de>0,$ then the solution $\tilde M_{\tilde t}$ of the
rescaled equation converges for $\tilde t\to \infty$ to a limiting submanifold $\tilde M_\infty$
satisfying the equation
$$F^\perp =H.$$
\end{thm}

\bibliographystyle{amsplain}

\end{document}